\documentclass[12pt]{amsart}
\usepackage{amsfonts}

\newtheorem{theorem}{Theorem}
\newtheorem{lemma}[theorem]{Lemma}
\newtheorem{proposition}[theorem]{Proposition}
\newtheorem{corollary}[theorem]{Corollary}

\theoremstyle{definition}
\newtheorem{definition}[theorem]{Definition}
\newtheorem{example}[theorem]{Example}
\newtheorem{notrems}[theorem]{Notation and Remarks}
\newtheorem{notation}[theorem]{Notation}
\newtheorem{remark}[theorem]{Remark}

\newcommand{\Section}[1]{\section{#1}\setcounter{theorem}{0}}

\newcommand{\<}{\langle}
\renewcommand{\>}{\rangle}
\newcommand{\la}{\langle}
\newcommand{\ra}{\rangle}

\newcommand{\z}{{\mathfrak z}}

\renewcommand{\v}{{\mathfrak v}}

\newcommand{\n}{{\mathfrak n}}
\newcommand{\x}{{\mathfrak x}}
\newcommand{\y}{{\mathfrak y}}
\newcommand{\so}{{\mathfrak{so}}}
\newcommand{\G}{{\mathcal G}}
\newcommand{\M}{{\mathcal M}}
\renewcommand{\Gamma}{{\varGamma}}
\renewcommand{\L}{{\mathcal L}}
\newcommand{\R}{{\mathbb R}}
\newcommand{\C}{{\mathbb C}}
\newcommand{\Z}{{\mathbb Z}}
\newcommand{\Q}{{\mathbb Q}}
\newcommand{\N}{{\mathbb N}}
\renewcommand{\epsilon}{{\varepsilon}}
\renewcommand{\phi}{{\varphi}}
\renewcommand{\theta}{{\vartheta}}

\newcommand{\scp}{{\la\,\,,\,\ra}}
\newcommand{\liebr}{{[\,\,,\,]}}

\newcommand{\ad}{{\operatorname{ad}}}

\newcommand{\Id}{{\operatorname{Id}}}
\newcommand{\kernn}{{\operatorname{ker}}}

\renewcommand{\exp}{{\operatorname{exp}}}

\newcommand{\grad}{{\operatorname{grad}}}
\newcommand{\proj}{{\operatorname{proj}}}
\newcommand{\dimm}{{\operatorname{dim}}}

\newcommand{\spann}{{\operatorname{span}}}
\newcommand{\minzero}{{\setminus\{0\}}}

\newcommand{\inv}{^{-1}}
\newcommand{\restr}[1]{\lower0.4ex\hbox{$|$}\lower0.7ex
  \hbox{$\scriptstyle{#1}$}}

\begin{document}

\title[Integrability of geodesic flows and isospectrality]
{Integrability of geodesic flows and isospectrality of Riemannian manifolds}

\author{Dorothee Schueth}

\address{Institut f\"ur Mathematik, Humboldt-Universit\"at zu
Berlin, D-10099 Berlin, Germany}

\email{schueth@math.hu-berlin.de}


\thanks{The author was partially supported by DFG Sonderforschungsbereich~647.}

\begin{abstract}
We construct a pair of compact, eight-dimensional, two-step
Riemannian nilmanifolds $M$ and $M'$ which are isospectral for the
Laplace operator on functions and such that $M$ has completely
integrable geodesic flow in the sense of Liouville, while $M'$ has
not. Moreover, for both manifolds we analyze the structure of the
submanifolds of the unit tangent bundle given by to
maximal continuous families of closed geodesics with generic
velocity fields. The structure of these submanifolds turns out to
reflect the above (non)integrability properties. On the other
hand, their dimension is larger than that of the Lagrangian tori
in~$M$, indicating a degeneracy which might explain
the fact that the wave invariants do not distinguish an integrable
from a nonintegrable system here.
Finally, we show that for~$M$, the invariant
eight-dimensional tori which are foliated
by closed geodesics are dense in the unit tangent bundle, and
that both $M$ and $M'$ satisfy the so-called Clean Intersection
Hypothesis.
\end{abstract}

\maketitle

\Section{Introduction}
\label{sec:intro}

\noindent
The spectrum of a compact Riemannian manifold is defined as the
collection of eigenvalues of the Laplace operator acting on
functions, counted with multiplicities. Two manifolds are called
\emph{isospectral} if their spectra are equal.

Inverse spectral geometry deals with the question of how much information
the spectrum of a manifold provides about its geometry.
Classical tools for extracting geometrical information from the spectrum
are asymptotic expansions of the singularities of the heat trace
or the wave trace. The so-called heat invariants determine the
dimension, the volume, the total scalar curvature and a series of
other integrals  depending on the curvature tensor
of the manifold. The singularities
of the wave trace are contained in the set of lengths of closed
geodesics on the manifold. Asymptotic expansions of these
singularities near such a length yield, under suitable nondegeneracy
assumptions, geometric information on the set of closed geodesics
of this length; see the fundamental article by Duistermaat and
Guillemin~\cite{DG:1975} or, for example,
S.~Zelditch's article~\cite{Zelditch:1998} for more detailed
results under stronger assumptions.

Closed geodesics thus being at the focus of the wave invariants,
it is natural to ask to which extent integrability 
properties of the dynamical system
given by the geodesic
flow of a Riemannian manifold
are determined by spectral data. In the present article, we will,
more precisely, answer the question whether complete integrability
of the geodesic flow in the sense of Liouville is a property determined
by the Laplace spectrum on functions.

A Riemannian manifold~$M$ is said to have \emph{completely integrable geodesic
flow} in the sense of Liouville if there exist $n=\dimm\,M$ first integrals
$f_1,\ldots f_n\in C^\infty(T^*M)$ of the geodesic flow which Poisson commute
and are functionally independent. These notions are defined as follows:
Let $\omega$ be the canonical symplectic form on $T^*M$. For any smooth function~$f$
on~$T^*M$, define the associated Hamiltonian vectorfield~$X_f$ by
$df=\omega(\,.\,,X_f)$. Then two functions $f,h\in C^\infty(T^*M)$
\emph{Poisson commute} if $\omega(X_f,X_h)=0$, or equivalently,
if $f$ is constant along the integral curves of $X_h$.
A set $\{f_1,\ldots,f_n\}$ of smooth functions on $T^*M$ is called
\emph{functionally independent} if $df_1\wedge\ldots\wedge df_n\ne0$
on an open dense subset of~$T^*M$.

\setcounter{section}{3}
\setcounter{theorem}{9}
Our main result is:

\begin{theorem}
There exists a pair of compact closed isospectral Riemannian manifolds $M, M'$
such that $M$ has completely integrable geodesic flow, while $M'$ does not
have completely integrable geodesic flow.
\end{theorem}

\setcounter{section}{1}
\setcounter{theorem}{1}
In recent years, L.~Butler has studied the question of complete integrability
of geodesic flows on compact Riemannian nilmanifolds. These are manifolds
of the form $(\Gamma\backslash N,g)$, where $N$ is a simply connected nilpotent
Lie group, $\Gamma$ is a cocompact discrete subgroup of~$N$, and $g$ is a
Riemannian metric induced by a left invariant metric on~$N$.
In particular, Butler has established sufficient criteria both for complete
integrability and lack of complete integrability for the case of two-step
Riemannian nilmanifolds~\cite{Butler:2003a}, \cite{Butler:2003b}.
Applying his results, together with a certain construction
by C.~Gordon and E.~Wilson~\cite{GW:1997} of isospectral nilmanifolds,
we find a pair of compact eight-dimensional
two-step Riemannian nilmanifolds $M, M'$ with the properties claimed
in the above theorem.

As an aside, we will also see that there exist continuous isospectral
families of compact two-step Riemannian nilmanifolds each of which
has completely integrable geodesic flow; see Corollary~\ref{cor:isodefo}.

For the manifold~$M$ we will, apart from just applying Butler's integrability
condition, explicitly present a set of eight Poisson
commuting and functionally independent
first integrals for the geodesic flow (see Lemma~\ref{lem:firstints}),
following the lines of Butler's proof
of his criterion in~\cite{Butler:2003a}. Note that these first integrals
are $C^\infty$ but not analytic; in fact, complete integrability by commuting
analytic first integrals is impossible for compact nilmanifolds which are not tori,
due to a theorem by I.~Taimanov~\cite{Taimanov:1987}.

Moreover, for any closed unit speed geodesic whose velocity field satisfies
a certain genericity condition, we will analyze the structure
of the maximal continuous family
of closed orbits of the geodesic flow in the unit tangent bundle of~$M$, resp.~$M'$,
containing the orbit associated with the given geodesic; see Corollary~\ref{cor:struct}.
While the dimension of $M$ and $M'$ is eight, these subsets of the unit
tangent bundles turn out to be nine-dimensional submanifolds.
``Generically'' one would expect, at least in the case of~$M$, eight
dimensional submanifolds instead, namely, fibers of the (singular) $T^8$ foliation
of the unit tangent bundle whose existence is guaranteed by complete
integrability of the geodesic flow. The fact that the dimension of the
submanifolds foliated by generic closed geodesics is higher than expected provides
some explanation for why the wave invariants can, in this example, not
distinguish between a completely integrable geodesic flow and a nonintegrable one,
in the sense that some degeneracy is present here.

On the other hand, the structure of these submanifolds reflects in a nice way
the (non)integrability properties of the geodesic flow: While in the case
of~$M$, the submanifolds under consideration
can be decomposed into a one-parameter family
of invariant $T^8$ fibers which are level sets of the first integrals
and two singular $T^7$ fibers, the picture is topologically quite different
in the case of~$M'$: Here, the corresponding submanifolds
cannot be decomposed into invariant $T^8$ fibers, but instead (again up to two singular
fibers) into a one-parameter family of invariant fibers diffeomorphic to
$H^3\times T^5$, where $H^3$ is a compact three-dimensional
Heisenberg manifold (a certain
two-step nilmanifold not diffeomorphic to a torus).

We will also show that for both $M$ and $M'$, the set of initial
vectors of unit speed closed geodesics
satisfying the considered genericity condition is dense in the unit
tangent bundle. In particular, the invariant $T^8$ fibers in the unit
tangent bundle of~$M$ which are level sets of the first integrals
and are fibered by closed orbits are dense in the unit tangent bundle;
see Corollary~\ref{cor:dense}.

Finally, using a result by Ruth Gornet~\cite{Gornet:2005},
we show in Remark~\ref{rem:CIH}
that the so-called Clean Intersection Hypothesis
(a hypothesis always needed in order to establish a wave trace formula~\cite{DG:1975}),
is satisfied for both $M$ and~$M'$. So the reason for why the wave
trace fails to distinguish between the different integrability properties
of the geodesic flows of~$M$ and~$M'$ does not lie in a failure of
the Clean Intersection Hypothesis, but must be subtler.

In fact, the main reason seems to lie in the phenomenon already
mentioned above, namely, that the closed geodesics fill out larger dimensional
families than Lagrangian tori.
In section~8 of his excellent survey
paper~\cite{Zelditch:2004}, Zelditch explains why the
the most natural setting for
the ``Can one hear integrability?''
problem or any
positive conjecture in this respect is that of
``simple clean length spectrum'', where, in the context of a manifold with
completely integrable geodesic flow, ``simple'' means that for any given
length, the family of closed geodesics of that length forms a Lagrangian
submanifold which should, moreover, be connected. The natural
approach for showing
that any Riemannian manifold isospectral to a manifold with these properties and
completely integrable geodesic flow should
again have integrable geodesic flow --- this approach has indeed yielded results
in dimension two, see~\cite{Zelditch:2004} --- is the following:
The wave invariants determine the dimension of the submanifolds foliated by
closed geodesics, and they determine certain properties of the geodesic flow on these
submanifolds. The hope is that if these submanifolds are Lagrangian tori
in the reference manifold, then the wave invariants would force them
to be tori also in the second manifold, which in turn might be used to derive
integrability of its geodesic flow.

However, as we have seen, the
Lagrangian condition is certainly violated in our examples because the dimension of
the submanifolds in question is larger than half the dimension of the cotangent bundle.
Let us mention here, without discussing it in the paper, that the connectedness
condition is not satisfied either: For a given length, there will in general
be more than
one continuous family of closed geodesics of that length.

The paper is organized as follows:

In Sections~\ref{sec:isosp} and~\ref{sec:int}, we set out the necessary
framework on two-step Riemannian nilmanifolds and present the results
of Gordon and Wilson~\cite{GW:1997} and Butler~\cite{Butler:2003a},
\cite{Butler:2003b}, respectively, which we need for our construction.
In Example~\ref{exple} we define the pair of Riemanian manifolds $M$ and $M'$
whose isospectrality is established in Section~\ref{sec:isosp}, while
their integrability vs.~nonintegrability properties are proven in
Section~\ref{sec:int}.
In Section~\ref{sec:explic} we explicitly establish eight Poisson commuting,
functionally independent first integrals of the geodesic flow of~$M$.
In Section~\ref{sec:struct} we analyze the geometry of the submanifolds
of the unit tangent bundles of $M$ and $M'$ which are foliated by maximal
continuous families of closed geodesics
satisfying a certain genericity condition. We conclude by proving
the results mentioned above on the
density of closed geodesics and on the Clean Intersection
Hypothesis.

The author would like to thank Steve Zelditch for raising her attention to the question
of whether complete integrability of the geodesic flow is spectrally determined,
and for several inspiring conversations.

\Section{Isospectral two-step nilmanifolds}
\label{sec:isosp}

\noindent
Let $\v$ and $\z$ be euclidean vector spaces, each endowed with a fixed
inner product.

\begin{definition}
\label{def:twostep}\
Given the above data, one associates with any linear map
$j:\z\to\so(\v)$ the following:
\begin{itemize}
\item[(i)]
The two-step nilpotent metric Lie algebra $\n(j)$ with underlying
vector space $\v\oplus\z$, whose inner product is given by letting
$\v$ and $\z$ be orthogonal and taking the
given inner product on
each factor, and whose Lie bracket $\liebr^j$ is defined by
letting $\z$ be central, $[\v,\v]^j\subseteq\z$ and $\la
j(Z)X,Y\ra=\la Z,[X,Y]^j\ra$ for all $X,Y\in\v$ and $Z\in\z$.
\item[(ii)]
The two-step simply connected nilpotent Lie group $N(j)$ whose Lie algebra is
$\n(j)$, and the left invariant Riemannian metric $g(j)$ on $N(j)$ which
coincides with the chosen inner product on $\n(j)=T_e N(j)$.
\end{itemize}
\end{definition}

\begin{notrems}
\label{notrems:latt}\
\begin{itemize}
\item[(i)]
Note that the Lie group exponential map $\exp^j:\n(j)\to N(j)$ is
a diffeomorphism because $N(j)$ is simply connected and nilpotent.
\item[(ii)]
Since $N(j)$ is two-step nilpotent, the Campbell-Baker-Hausdorff formula
implies that $\exp^j(X)\cdot\exp^j(Y)=\exp^j(X+Y+\frac12[X,Y]^j)$ for all
$X,Y\in\n(j)$.
\item[(iii)]
In particular, if a lattice $\G$ in $\v\oplus\z$ has the property
that $[\G,\G]^j\subseteq 2\G$, then $\Gamma:=\exp^j(\G)$ is a discrete
subgroup of $N(j)$; hence $\Gamma\backslash N(j)$, endowed with the metric
induced by $g(j)$, is a two-step
Riemannian nilmanifold. We denote the induced metric by $g(j)$ again.
If, moreover, $\G$ has full rank, then $\Gamma$
is cocompact, and $(\Gamma\backslash N(j),g(j))$ is a compact two-step Riemannian
nilmanifold.
\item[(iv)]
If $\L$ is a cocompact lattice in~$\z$, then we denote by $\L^*:=\{Z\in\z\mid
\la Z,\L \ra\subseteq\Z$\} the dual lattice, viewed as a lattice in $\z$.
\end{itemize}
\end{notrems}

\begin{definition}
\label{def:jLisosp}\
\begin{itemize}
\item[(i)]
Two linear maps $j,j':\z\to\so(\v)$ are called \emph{isospectral}
if for each $Z\in\z$, the maps $j(Z),j'(Z)\in\so(\v)$ are similar,
that is, have the same eigenvalues (with multiplicities) in~$\C$.
\item[(ii)]
Two lattices in a euclidean vector space are called \emph{isospectral}
if the lengths of their elements, counted with multiplicities,
coincide.
\end{itemize}
\end{definition}

We will use the following special version of a result from \cite{GW:1997}.

\begin{proposition}[see \cite{GW:1997} 3.2, 3.7, 3.8]
\label{prop:gw}
Let $j,j':\z\to\so(\v)$ be isospectral. Let $\M$ and~$\L$ be cocompact lattices
in $\v$ and~$\z$, respectively. Assume that $[\M,\M]^j$ and $[\M,\M]^{j'}$
are contained in $2\L$. For each $Z\in\L^*$ assume that the lattices $\kernn
(j(Z))\cap\M$ and $\kernn(j'(Z))\cap \M$ are isospectral.
Write $\Gamma(j):=\exp^j(\M+\L)$, and similarly for $j'$.
Then the compact Riemannian manifolds $(\Gamma(j)\backslash N(j),g(j))$
and $(\Gamma(j')\backslash N(j'),g(j'))$ are isospectral for the Laplace
operator on functions.
\end{proposition}

\begin{remark}
\label{rem:proofgw}
(i)
In the situation of Proposition \ref{prop:gw}, note that for $Z\in\L^*$, the
subspace $\la Z\ra^\perp$ of $\z$ is rational with respect to~$\L$, that is, its intersection
with $\L$ is a cocompact lattice in this subspace. Using the assumption
$[\M,\M]^j\subseteq 2\L\subset\L$ and the fact that $\M$ has
full rank in~$\v$, one concludes that the subspace
$\{X\in\v\mid[X,\v]^j\subseteq\la Z\ra^\perp\}=\kernn(j(Z))$ of~$\v$ is rational
with respect to $\M$. Thus
$\kernn(j(Z))\cap\M$ is actually a cocompact lattice in $\kernn(j(Z))$, and
similarly for~$j'$.

(ii)
We indicate how to derive Proposition \ref{prop:gw} from the cited results of \cite{GW:1997}.

Let $Z\in\L^*\minzero$. Consider the compact two-step Riemannian
nilmanifold $M_Z$ associated as in \ref{def:twostep} and
\ref{notrems:latt} with the euclidean spaces $\v$ and $\la Z\ra$,
the linear map from $\la Z\ra$ to $\so(\v)$ mapping $Z$ to $j(Z)$,
and the lattices $\M\subset\v$ and $\proj_{\la Z\ra}\L\subset \la
Z\ra$, where $\proj_{\la Z\ra}$ denotes orthogonal projection
onto~$\la Z\ra$. Analogously define $M'_Z$ using $j(Z)$ instead of
$j'(Z)$. ($M_Z$~is actually the Riemannian submersion quotient of
$(\Gamma(j)\backslash N(j),g(j))$ by the canonical action of the
torus $\la Z\ra^\perp/(\la Z\ra^\perp\cap\L)$
on this manifold, and
similarly for $j'$.)

We first note that $M_Z$ and $M'_Z$ are then isospectral by
Proposition~3.7/Remark~3.8 of~\cite{GW:1997}. In fact, their base tori
(the Riemannian submersion quotients by the action of the circle $\la Z\ra/
(\la Z\ra\cap\proj_{\la Z\ra}\L)$) are both isometric to
the torus $\v/\M$ and thus isospectral, the maps $j(Z)$ and $j'(Z)$
are similar, the lattice in $\kernn(j^{(\prime)}(Z))\oplus\la Z\ra$
is a sum of lattices in the two factors as required in Remark~3.8 of~\cite{GW:1997},
and the lattices $\ker(j(Z))\cap\M$ and $\ker(j'(Z))\cap\M$ are isospectral by assumption.

Finally, the manifold $M_0$, defined as the Riemannian submersion quotient
of the manifold
$(\Gamma(j)\backslash N(j),g(j))$ by the action of the torus $\z/\L$,
is isometric to the analogously defined manifold~$M'_0$
since both are isometric to $\v/\M$;
in particular, $M_0$ and $M'_0$ are isospectral.
Isospectrality of the pairs $M_Z$ and $M'_Z$ for each $Z\in\L^*$ now implies,
by Theorem~3.2 of~\cite{GW:1997},
the isospectrality statement of Proposition~\ref{prop:gw}.
\end{remark}

We now give an example of a pair of isospectral manifolds
arising from Proposition~\ref{prop:gw} with the property that, as we will
show in Section~\ref{sec:int}, the geodesic flow on the
first manifold is completely integrable,
while the geodesic flow on the second manifold is not.

\begin{example}
\label{exple} Let $\dimm\,\v=5$, $\dimm\,\z=3$, and let
$\{X_i,X_j,Y_i,Y_j,Y_k\}$ and $\{Z_i,Z_j,Z_k\}$ be orthonormal
bases of $\v$ and~$\z$, respectively. Define two Lie brackets
$\liebr$ and~$\liebr'$ on $\v\oplus\z$ as follows: For
$a,b\in\{i,j,k\}$ with $a\ne b$ denote by $ab\in\{\pm i,\pm j,\pm
k\}$ the quaternionian product of $a$ and~$b$, and for
$c\in\{i,j,k\}$ write $Z_{-c}:=-Z_c$. Now let
$[X_a,Y_b]=-[Y_b,X_a]:=Z_{ab}$ and $[X_a,X_b]':=Z_{ab}$,
$[Y_a,Y_b]':=Z_{ab}$ for all $a\ne b$ in $\{i,j,k\}$, ignoring any
expressions containing $X_k$, and let all other brackets between
basis elements be zero. Then $\liebr=\liebr^j$ and
$\liebr'=\liebr^{j'}$, where $j,j':\z\to\so(\v)$ are the linear
maps for which $j(c_iZ_i+c_jZ_j+c_kZ_k)$ and $j'(c_iZ_i+
c_jZ_j+c_kZ_k)$ with $c_i,c_j,c_k\in\R$ are expressed, with
respect to the given basis of~$\v$, by the matrices
$$
\left(\begin{smallmatrix}
0&0&0&-c_k&c_j\\ 0&0&c_k&0&-c_i\\0&-c_k&0&0&0\\c_k&0&0&0&0\\-c_j&c_i&0&0&0
\end{smallmatrix}\right)\quad\mathrm{and}\quad
\left(\begin{smallmatrix}
0&-c_k&0&0&0\\c_k&0&0&0&0\\0&0&0&-c_k&c_j\\0&0&c_k&0&-c_i\\0&0&-c_j&c_i&0
\end{smallmatrix}\right),\quad\mathrm{respectively.}
$$
Let $\M$ be the lattice in~$\v$ generated by the basis given above,
and $\L$ be the lattice in~$\z$ generated by $\{\frac12Z_i,\frac12Z_j,\frac12Z_k\}$;
in particular, $[\M,\M]^j=[\M,\M]^{j'}=2\L$.
Then, with notation as in \ref{def:twostep}, \ref{notrems:latt}, \ref{prop:gw},
the associated manifolds $(\Gamma(j)\backslash N(j),g(j))$ and
$(\Gamma(j')\backslash N(j'),g(j'))$ are isospectral by Proposition~\ref{prop:gw}.

In fact, $j$ and $j'$ are isospectral since the two matrices given above
have the same characteristic polynomials $\lambda(\lambda+c_k^2)(\lambda+c_i^2
+c_j^2+c_k^2)$. It remains to show that for each $Z\in\L^*$, the two lattices
$\kernn(j(Z))\cap\M$ and $\kernn(j'(Z))\cap\M$ are isospectral; we will actually
see that they are isometric.
For $c=(c_i,c_j,c_k)\in\R^3$ write $X_c:=c_iX_i+c_jX_j$,
$Y_c:=c_iY_i+c_jY_j+c_kY_k$, $Z_c:=c_iZ_i+c_jZ_j+c_kZ_k$.
Then we have for all $c\in\R^3$ and, in particular, for all $c\in(2\Z)^3$
(that is, for $Z_c\in\L^*$):
\begin{equation*}
\begin{split}
\kernn(j(Z_c))&=\kernn(j'(Z_c))=\spann\{Y_c\}\quad\mathrm{if\ }c_k\ne0,\\
\kernn(j(Z_c))&=\spann\{X_c,Y_i,Y_j\},\quad\kernn(j'(Z_c))=\spann\{X_i,X_j,Y_c\}
\quad\mathrm{if\ }c_k=0\mathrm{\ and\ }c\ne0,\\
\kernn(j(Z_c))&=\kernn(j'(Z_c))=\v\quad\mathrm{if\ }c=0.
\end{split}
\end{equation*}
In the first case, the two subspaces are equal, and thus the two lattices
in question coincide. The same holds for the third case.
In the second case, note that $Y_c=c_iY_i+c_jY_j$ since
$c_k=0$; this and the fact that $\M\cap\spann\{X_i,X_j,Y_i,Y_j\}$ is invariant
under exchanging the $X$- with the $Y$-space shows isometry of the lattices
in the two subspaces also in this case.
\end{example}

\begin{remark}
The metric Lie algebras $\n$ and $\n'$ are not isomorphic.
This implies~\cite{GW:1997} that the two manifolds
$(\Gamma\backslash N,g)$
and $(\Gamma'\backslash N',g')$ from Example~\ref{exple}
are not locally isometric.
\end{remark}

\Section{Complete integrability of geodesic flows on two-step nilmanifolds}
\label{sec:int}

\noindent
We are going to use results by L.~Butler \cite{Butler:2003a}, \cite{Butler:2003b}
in order to show that the first of the two isospectral manifolds from
Example~\ref{exple} has completely integrable geodesic flow, while
the second has not.

\begin{definition}
\label{def:hr}
\begin{itemize}
\item[(i)]
A two-step nilpotent Lie algebra $\n$ is called a \emph{Heisenberg-Reiter
Lie algebra}, shortly: \emph{HR Lie algebra}, if there exists a vector
space decomposition $\n=\x\oplus\y\oplus\z$
such that $[\n,\n]\subseteq\z$, $[\z,\n]=0$,
$[\x,\x]=0$, and $[\y,\y]=0$.
Such a decomposition is then called a \emph{presentation} of~$\n$.
\item[(ii)]
A presentation $\x\oplus\y\oplus\z$ of an HR Lie algebra~$\n$
is called \emph{injective} if there exists $c\in\z^*$ such that
$c\restr{[X,\y]}\ne0$ for all $X\in\x\minzero$.
\end{itemize}
\end{definition}

\begin{theorem}[\cite{Butler:2003a}, Theorem~2.22]
\label{butlerint}
Let $\n$ be an HR Lie algebra admitting an injective presentation,
and let $N$ be the associated simply connected Lie group.
Assume that there exists a discrete, cocompact subgroup $\Gamma$ of~$N$.
Then for any such $\Gamma$
and any left invariant metric~$g$ on~$N$, the geodesic flow of
$(\Gamma\backslash N,g)$ is completely integrable.
\end{theorem}

\begin{corollary}
\label{cor:int}
The manifold $(\Gamma(j)\backslash N(j),g(j))$ from Example~\ref{exple}
has completely integrable geodesic flow.
\end{corollary}

\begin{proof}
We use the notation from Example~\ref{exple}.

Let $\x:=\spann\{X_i,X_j\}$ and $\y:=\spann\{Y_i,Y_j,Y_k\}$.
Since $\x$ and $\y$ are abelian, $\n(j)=\x\oplus\y\oplus\z$ is
an HR Lie algebra. Note that $j(Z_k)\restr{\x}$ is injective, and
hence $\la Z_k,[X,\y] \ra = \la j(Z_k)X,\y \ra\ne0$ for all $X\in\x\minzero$.
It follows that $\x\oplus\y\oplus\z$ is an injective representation
of~$\n(j)$, and thus the corollary follows from Theorem~\ref{butlerint}.
\end{proof}

\begin{remark}
\label{rem:problem}
\begin{itemize}
\item[(i)]
It seems that the proof of Theorem~\ref{butlerint}
in~\cite{Butler:2003a} presupposes
that the rational and the HR structure of the Lie algebra~$\n$
are in some sense compatible; more precisely: The proof given there
of complete
integrability of the geodesic flow of $(\Gamma\backslash N,g)$
works without any additional arguments if there exists an injective
presentation $\x\oplus\y\oplus\z$ of~$\n$ with the property
that there exists a complement $\tilde\x$ of $\z$ in $\x\oplus\z$
and a complement $\tilde\y$ of $\z$ in $\y\oplus\z$ such that
the union $\tilde\x\cup\tilde\y\cup\z$ contains a set
of vectors which is mapped to a generating set for~$\Gamma$
by the Lie group exponential map. (Of course, $\tilde\x\oplus
\tilde\y\oplus\z$ is then itself an injective presentation of~$\n$.)
\item[(ii)]
By our construction of $\Gamma(j)$ in Example~\ref{exple} it is
clear that the injective representation $\x\oplus\y\oplus\z$ of
$\n(j)$ from the proof of Corollary~\ref{cor:int} does satisfy the
additional assumption just mentioned. Moreover, we will reconfirm
complete integrability of the geodesic flow on
$(\Gamma(j)\backslash N(j),g(j))$ in Section~\ref{sec:explic}
directly: Guided by the proof of Theorem~\ref{butlerint}
in~\cite{Butler:2003a}, we will explicitly establish eight
commuting, functionally independent first integrals.
\end{itemize}
\end{remark}

\begin{remark}
\label{rem:isodefo}
It is easy to find examples of continuous families of isospectral
manifolds each of which has completely integrable geodesic flow,
using Butler's above result and not the above construction from~\cite{GW:1997},
but another construction by Gordon and Wilson~\cite{GW:1984}, involving so-called
almost-inner automorphisms. For example, let $\dimm\,\v=4$, $\dimm\,\z=2$,
and let $\{X_1,X_2,Y_1,Y_2\}$ and $\{Z_1,Z_2\}$ be orthonormal bases of
$\v$ and~$\z$, respectively. Define a Lie bracket on the orthogonal sum
$\n:=\v\oplus\z$ by letting
$[X_1,Y_1]=[X_2,Y_2]=Z_1$, $[X_1,Y_2]=Z_2$, and letting pairs of basis vectors
commute if they do not occur in these three equations. 
Let $N$ be the associated simply connected Lie group, associated with the
left invariant metric~$g$ defined by the given inner product on~$\n$.
Now consider the continuous family of lattices $\Gamma_t:=\exp(\G_t)$ in the
associated simply connected Lie group~$N$,
where $\G_t:=\spann\{X_1,X_2,Y_1,Y_2+tZ_2,\frac12Z_1,\frac12Z_2\}$.
By~\cite{GW:1984}, the family $(\Gamma_t\backslash N,g)$ is nontrivial
and isospectral
(even strongly isospectral, that is, also on $p$-forms for all~$p$).
Note that these manifolds are locally isometric to each other
because it is just the lattice that changes. Using the
automorphism~$\Phi_t$ of~$N$ whose differential on~$\n$ is
given by $\Id+\<\,.\,,Y_2\>\cdot tZ_2$, we can instead
view this family
as the deformation of the metric on a fixed underlying manifold
because $(\Gamma\backslash N,\Phi_t^*g)$ is isometric to
$(\Gamma_t\backslash N,g)$.

Since
$\x:=\spann\{X_1,X_2\}$ and $\y:=\spann\{Y_1,Y_2\}$ are abelian, $\n$~is
an HR Lie algebra. Moreover, $j(Z_1)$ is invertible; in particular,
$\x\oplus\y\oplus\z$ is an injective presentation of~$\n$.
Letting $\tilde\x:=\x$ and $\tilde\y:=\spann\{Y_1,Y_2+tZ_2\}$ in
Remark~\ref{rem:problem}(ii),
we see that Butler's integrability theorem applies; hence, each of the manifolds
in this isospectral family has completely integrable geodesic flow.
It is also an easy exercise to find six commuting first integrals for each
of these manifolds, along the lines of Butler's proof, similarly to how
we will do this in the next
section for the first manifold from Example~\ref{exple}.
So we can state:
\end{remark}

\begin{corollary}
\label{cor:isodefo}
There exist continuous isospectral families of compact Riemannian
manifolds each of which has completely integrable geodesic flow.
\end{corollary}

However, it remains an open question whether there might even
exist a continuous isospectral
deformation of a manifold with completely integrable geodesic flow
to a manifold whose geodesic flow is nonintegrable.
Note that in the family just constructed, the geodesic flows of
the manifolds involved, even though they share the property of
complete integrability,
are not symplectically conjugate (i.e., conjugate by a symplectomorphism)
to each other because any pair of compact two-step Riemannian nilmanifolds
with symplectically conjugate geodesic flows must be
isometric~\cite{GSM:1997}.

We now consider a sufficient nonintegrability criterion by Butler~\cite{Butler:2003b}:

\begin{definition}
\label{def:noncomm}
Let $\n$ be a two-step nilpotent Lie algebra.
\begin{itemize}
\item[(i)]
For $\lambda\in\n^*$ let
$\n_\lambda:=\{X\in\n\mid\ad_X^*\lambda=0\}
=\{X\in\n\mid\lambda\restr{[X,\n]}=0\}$.
\item[(ii)]
$\lambda\in\n^*$ is called \emph{regular} if $\n_\lambda$ has
minimal dimension.
\item[(iii)]
$\n$ is called \emph{nonintegrable}
if there exists a dense open subset $U$ of $\n^*\times\n^*$ such that
for each $(\lambda,\mu)\in U$, both $\lambda$ and $\mu$ are
regular and $[\n_\lambda,\n_\mu]$ has positive dimension.
\end{itemize}
\end{definition}

\begin{theorem}[\cite{Butler:2003b}, Theorem 1.3]
\label{butlernonint}
Let $\n$ be a nonintegrable two-step nilpotent Lie algebra,
and let $N$ be the associated simply connected Lie group.
Assume that there exists a discrete, cocompact subgroup $\Gamma$ of~$N$.
Then for any such~$\Gamma$
and any left invariant metric~$g$ on~$N$, the geodesic flow of
$(\Gamma\backslash N,g)$ is \emph{not} completely integrable.
\end{theorem}

\begin{corollary}
\label{cor:nonint}
The manifold $(\Gamma(j')\backslash N(j'),g(j'))$ from
Example~$\mathrm{\ref{exple}}$
does not have completely integrable geodesic flow.
\end{corollary}

\begin{proof}
We use the notation from Example~\ref{exple}. Letting
$\n':=\n(j')$, we write elements
of $\n'=\v\oplus\z$ in the form $V+Z$ with $V\in\v$ and $Z\in\z$.
Then for $(V+Z)^*:=\la V+Z,\,.\,\ra\in\n^{\prime\,*}$ we have
$\n'_{(V+Z)^*} = \kernn(j'(Z))\oplus\z$.
Thus $(V+Z)^*$ is regular if and only if $Z=Z_c$ for some $c\in\R^3$ with $c_k\ne0$;
in that case, $\n'_{(V+Z)^*} = \spann\{Y_c\}\oplus\z$
(see the discussion in Example~\ref{exple}).
Moreover, if $\tilde V\in\v$ and $\tilde Z=Z_{\tilde c}$ with $\tilde c_k\ne0$, then
$$
[\n'_{(V+Z)^*},\n'_{(\tilde V+\tilde Z)^*}]=\spann\{[Y_c\,,Y_{\tilde c}]^{j'}\},
$$
which has positive dimension (namely, dimension one) if and only if $c$ and $\tilde c$
are linearly independent, or equivalently: if $Z$ and $\tilde Z$ are linearly independent.
The set of pairs of vectors $(V+Z,\tilde V+\tilde Z)\in\n'\times\n'$
such that $Z$ and $\tilde Z$ are linearly independent and both have nonvanishing
$Z_k$-component is obviously open and dense in $\n'\times\n'$.
Using the identification of $\n'$ with $\n^{\prime\,*}$ induced by $\scp$, one concludes
that there is an open and dense subset $U$ of $\n^{\prime\,*}
\times\n^{\prime\,*}$
with the property required in Definition~\ref{def:noncomm}(iii). Thus $\n'$ is
nonintegrable, and the corollary follows from Theorem~\ref{butlernonint}.
\end{proof}

We now conclude our main result:

\begin{theorem}
There exists a pair of isospectral compact closed Riemannian manifolds
$M,M'$ such that $M$ has completely integrable geodesic flow,
while $M'$ does not have completely integrable geodesic flow.
\end{theorem}

\begin{proof}
This follows immediately from Example~\ref{exple} and
Corollaries~\ref{cor:int} and~\ref{cor:nonint},
letting $M:=(\Gamma(j)\backslash N(j),g(j))$
and $M':=(\Gamma(j')\backslash N(j'),g(j'))$.
\end{proof}

\Section{Explicit first integrals}
\label{sec:explic}

\noindent In this section we will explicitly establish
eight Poisson commuting first integrals for 
the geodesic flow of the first manifold $(\Gamma(j)\backslash
N(j),g(j))$ from Example~\ref{exple},
using the ideas of the proof of
Theorem~\ref{butlerint} from~\cite{Butler:2003a} (see also
Remark~\ref{rem:problem}). We will do our computations on the
tangent bundle rather than on the cotangent bundle; these two are
canonically identified by the euclidean metric induced by $g(j)$
on each tangent space.

\begin{notrems}
\label{not:twostep} Let $\v$, $\z$, $j$ be data as in
Definition~\ref{def:twostep}, let $\n:=\n(j)$ be the associated
two-step nilpotent Lie algebra with underlying vector space
$\v\oplus\z$ and Lie bracket $\liebr:=\liebr^j$, let $N:=N(j)$ be
the associated simply connected Lie group, $g:=g^j$ be the
associated left invariant metric on~$N$, and $\exp:=\exp^j:\n\to
N$ be the Lie group exponential map.
\begin{itemize}
\item[(i)]
We denote elements of $\n$ in the form $V+Z$ with $V\in\v$ and
$Z\in\z$, and we denote elements of $N$ in the form
$(v,z):=\exp(v+z)$ with $v\in\v$ and $z\in\z$. By the
Campbell-Baker-Hausdorff formula we have
$$(v,z)\cdot(\bar v,\bar z)=(v+\bar v,z+\bar
z+\textstyle{\frac12}[v,\bar v]).
$$
\item[(ii)]
For $a\in N$, we denote left multiplication by $a$
by $L_a:N\to N$. We identify the tangent
bundle $TN$ with $N\times\n$ using left translation; that is, we
write $X\in T_{(v,z)}N$ in the form $((v,z),L\inv_{(v,z)*}X)\in
N\times\n$. Note that for $a\in N$ we have
$$
L_{a*}((v,z),X)=(L_a(v,z),X).
$$
\end{itemize}
\end{notrems}

\begin{lemma} \cite{Eberlein:1994}
\label{lem:geodeqs} In the situation
of~$\,\mathrm{\ref{not:twostep}}$, let $\gamma:\R\to N$ be a
geodesic in $(N,g)$. Write $\dot\gamma(t)= \bigl((v(t),z(t)),
V(t)+Z(t)\bigr)$ with $v(t),V(t)\in\v$ and $z(t),Z(t)\in\z$. Then
the following \emph{geodesic equations} hold:
\begin{equation*}
\begin{gathered}
\dot V(t)=j(Z(t))V(t),\quad \dot Z(t)=0,\\
\dot v(t)=V(t),\quad \dot
z(t)=Z(t)+\textstyle{\frac12}[v(t),V(t)];
\end{gathered}
\end{equation*}
 hence $Z(t)\equiv Z(0)=:Z$ and
\begin{gather*}
V(t)=e^{tj(Z)}V(0),\\
\dot v(t)=e^{tj(Z)}V(0),\quad \dot z(t)=Z+\textstyle{\frac12}
[v(t),V(t)].
\end{gather*}
\end{lemma}

\begin{notrems}
\label{not:ourj}
In the following, let $\v$, $\z$ be as in
Example~\ref{exple}, and let $j:\z\to\so(\v)$ be the first of the
two maps considered there.
\begin{itemize}
\item[(i)]
For $Z\in\z$ with $Z=Z_c=c_iZ_i+c_jZ_j+c_kZ_k$ let
\begin{gather*}
E_1(Z):=c_iX_i+c_jX_j,\quad E_2(Z):=-c_jY_i+c_iY_j\\
E_3(Z):=|c|(c_jX_i-c_iX_j), \quad
E_4(Z):=c_k(c_iY_i+c_jY_j)-(c_i^2+c_j^2)Y_k
\end{gather*}
and
$$Y(Z):=Y_c=c_iY_i+c_jY_j+c_kY_k.
$$
Note that we then have
\begin{equation}\label{jeq}
\begin{gathered}
j(Z)E_1(Z)=c_kE_2(Z),\quad
j(Z)E_2(Z)=-c_kE_1(Z),\\
j(Z)E_3(Z)=|c|E_4(Z),\quad j(Z)E_4(Z)=-|c|E_3(Z),\\
j(Z)Y(Z)=0.
\end{gathered}
\end{equation}
Thus, for generic $Z$ (namely, with $|c|>|c_k|>0$), the sets
$\{E_1(Z),E_2(Z)\}$ and $\{E_3(Z),E_4(Z)\}$ are bases of the
eigenspaces associated with the eigenvalues $-c_k^2$ and $-|c|^2$
of $j(Z)^2$, respectively, and $Y(Z)$ spans the zero eigenspace of
$j(Z)$.
\item[(ii)]
Let $\x:=\spann\{X_i,X_j\}$, $\y:=\spann\{Y_i,Y_j,Y_k\}$ as in the
proof of Corollary~\ref{cor:int}. For $Z\in\z$ with
$Z=c_iZ_i+c_jZ_j+c_kZ_k$ and $c_k\ne0$ we denote by $C(Z):\y\to\x$
the linear map whose matrix with respect to the given bases of
$\y$ and $\x$ is
$$
(c_k|c|^2)\inv\left(\begin{smallmatrix}
-c_ic_j&c_i^2+c_k^2&-c_jc_k\\-c_j^2-c_k^2&c_ic_j&c_ic_k\end{smallmatrix}\right)
$$
Note that we then have
$$C(Z)\circ j(Z)\restr\x=\Id_\x\,.
$$
In fact, $C(Z)$ is just $(j(Z)\restr\y\,j(Z)\restr\x)\inv
j(Z)\restr\y$\,.
\item[(iii)]
Let $\phi\in C^\infty(\R)$ be the map $x\mapsto e^{-1/x^2}$ with
$\phi(0):=0$. Define $\Phi:\z\to\R$ by
$\Phi(c_iZ_i+c_jZ_j+c_kZ_k):=\phi(c_k|c|^2)$.
\item[(iv)]
For $V\in\v=\x\oplus\y$ denote by $V_\x$ and $V_\y$ the components
of~$V$ in $\x$ and~$\y$, respectively.
\end{itemize}
\end{notrems}

\begin{lemma}
\label{lem:firstints}
We use
Notation~$\,\mathrm{\ref{not:twostep}}$. Then in the situation
of~$\,\mathrm{\ref{not:ourj}}$, each of the following eight
functions on~$TN$ is a smooth first integral of the geodesic flow
on $(N,g):=(N(j),g(j))$ and is invariant under the left action of
$\Gamma:=\Gamma(j)$ from Example~$\,\mathrm{\ref{exple}}$:
\begin{align*}
q^W:TN\ni((v,z),V+Z)&\mapsto\<Z,W\>\in\R\quad\mathrm{with\ }W\in\{Z_i,Z_j,Z_k\},\\
h_1:TN\ni((v,z),V+Z)&\mapsto\<V,E_1(Z)\>^2+\<V,E_2(Z)\>^2\in\R,\\
h_2:TN\ni((v,z),V+Z)&\mapsto\<V,E_3(Z)\>^2+\<V,E_4(Z)\>^2\in\R,\\
k:TN\ni((v,z),V+Z)&\mapsto\<V,Y(Z)\>\in\R,\\
f^X:TN\ni((v,z),V+Z)&\mapsto
\begin{cases}0\in\R\mathrm{\ if\ }Z=Z_c\mathrm{\ with\ }c_k=0,\\ \Phi(Z)\sin(
2\pi\<X,v_\x-C(Z)V_\y\>)\in\R,\mathrm{\ else,}
\end{cases}\\
&\mathrm{with\ }X\in\{X_i,X_j\}.
\end{align*}
In particular, each of these functions descends to a first
integral of the geodesic flow on $(\Gamma\backslash N,g)$.
\end{lemma}

\begin{proof}
Smoothness of the above functions is immediate from their construction.
Note that $q^W$, $h_1$, $h_2$, and $k$ are invariant under the
left action of $N$ on~$TN$. Moreover, if $a=(\bar v,\bar
z)\in\Gamma$ then $\bar v$ is an integer combination of the basis
vectors; hence
$$f^X\bigl(L_{a*}((v,z),V+Z)\bigr)=f^X\bigl((\bar v+v,\bar
z+z+\textstyle{\frac12}[\bar v,v]),V+Z\bigr)=f^X\bigl((v,z),V+Z\bigr)
$$
for $X\in\{X_i,X_j\}$ because $\<X,v_\x\>$ differs from
$\<X,v_\x+\bar v_\x\>$ by an integer.

It remains to show that each of the eight functions is invariant
under the geodesic flow. Let $\gamma$ be a geodesic in $(N,g)$ and
write $\gamma(t)=\bigl((v(t),z(t)),V(t)+Z(t)\bigr)$. Then
$Z(t)=:Z$ is constant by Lemma~\ref{lem:geodeqs}; in particular,
$q^W\circ\dot\gamma$ is constant. Always using~\ref{lem:geodeqs}
and the equations~(\ref{jeq}), we observe:
\begin{align*}
(h_1\circ\dot\gamma)'(t)&=2\<V(t),E_1(Z)\>\<j(Z)V(t),E_1(Z)\>
+2\<V(t)E_2(Z)\>\<j(Z)V(t),E_2(t)\>\\ &=
-2\<V(t),E_1(Z)\>\<V(t),c_k
E_2(Z)\>-2\<V(t),E_2(t)\>\<V(t),-c_kE_1(Z)\>=0, \end{align*} and
similarly for $h_2$. Moreover,
$$
(k\circ\dot\gamma)'(t)=\<j(Z)V(t),Y(Z)\>=-\<V(t),j(Z)Y(Z)\>=0.
$$
Finally, noting that $\frac
d{dt}V(t)_\y=(j(Z)V(t))_\y=j(Z)V(t)_\x$\,, we have
$$
(f^X\circ\dot\gamma)'(t)=\Phi(Z)\cos(
2\pi\<X,v(t)_\x-C(Z)V(t)_\y)\cdot
2\pi\<X,V(t)_\x-C(Z)j(Z)V(t)_\x\>=0
$$
if $Z=Z_c$ with $c_k\ne0$; if $c_k=0$ then $(f^X\circ\dot\gamma)(t)\equiv0$
by definition.
\end{proof}

\begin{lemma}
\label{poiscom}
\begin{itemize}
\item[(i)] The eight first integrals from Lemma~$\mathrm{\,\ref{lem:firstints}(i)}$
are functionally independent, and
\item[(ii)] they Poisson commute with each other.
\end{itemize}
\end{lemma}

\begin{proof}
(i) Note that with respect to the left invariant Riemannian
product metric $g\times\scp$ on $TN\cong N\times\n$, the gradients
of the functions $q^{Z_i}, q^{Z_j}, q^{Z_k},h_1,h_2,k$ are all
tangent to the second factor. The gradients of $q^{Z_i}, q^{Z_j},
q^{Z_k}$ at $((v,z),V+Z)\in TN$ are just
\begin{equation*}
(0,Z_i), (0,Z_j), (0,Z_k)\in T_{(v,z)}N\oplus T_{V+Z}\n\cong T_{(v,z)}N\oplus\n.
\end{equation*}
The gradients of $h_1, h_2, k$ at $((v,z),V+Z)$, viewed as elements
of $T_{(v,z)}N\oplus\n$, are of the form
\begin{equation*}
\begin{gathered}
(0,\,2\<V,E_1(Z)\>E_1(Z)+2\<V,E_2(Z)\>E_2(Z)+W),\\
(0,\,2\<V,E_3(Z)\>E_3(Z)+2\<V,E_4(Z)\>E_4(Z)+\tilde W),\\
(0,Y(Z)+U),
\end{gathered}
\end{equation*}
respectively, where $W,\tilde W, U$ are in $\z$. If $Z=Z_c$ with
$|c|>|c_k|>0$, and if $V$ is not orthogonal to any of the
subspaces $\spann\{E_1(Z),E_2(Z)\}$, $\spann\{E_3(Z),E_4(Z)\}$,
and $\spann\{Y(Z)\}$, then these six gradients are obviously linearly
independent (recall~\ref{not:ourj}(i)). Moreover, for these~$Z$,
the gradient of $f^{X_i}$ at the point $((v,z),V+Z)$ is of the
form
\begin{equation*}
(0,W)+\Phi(Z)\cos(2\pi\<X_i,v_\x-C(Z)V_\y\>)
\cdot2\pi(L_{(v,z)*}X_i\,,\,-{}^t C(Z)X_i) \in T_{(v,z)}N\oplus\n
\end{equation*}
with some $W\in\z$, and similarly for $X_j$. Since $L_{(v,z)*}X_i$
and $L_{(v,z)*}X_j$ are linearly independent, it follows that the
set of points in $TN$ at which all eight gradients are linearly
independent is open and dense in~$TN$.

(ii) The symplectic form~$\omega$ on $TN$, after identification
with $T^*N$ by the left invariant metric~$g$, is given at the
point $((v,z),V+Z)\in TN\cong N\times\n$ by
\begin{align*}
\omega_{((v,z),V+Z)}\bigl((L_{(v,z)*}A,B),(L_{(v,z)*}\tilde
A,\tilde B)\bigr)&= \<B,\tilde A\>-\<A,\tilde B\>-\<V+Z,[A,\tilde
A]\>\\ &= \<B,\tilde A\>-\<A,\tilde B\>-\<j(Z)A_\v\,,\tilde A_\v\>
\end{align*}
for all $A,B,\tilde A,\tilde B\in\n$, where $A_\v$ denotes the
$\v$-component of $A$ in $\n=\v\oplus\z$. For any $f\in
C^\infty(TN)$, the Hamiltonian vectorfield $X_f$ is characterized
by $\<\,.\,,\grad f\>=df=\omega(\,.\,,X_f)$. This and the above
formula for~$\omega$ shows that if the gradient of~$f$ at the
point $((v,z),V+Z)\in TN$ is $(L_{(v,z)*}B,A)$, then $X_f$ at this
point is
$$(L_{(v,z)*}A,-B+j(Z)A_\v).
$$
Therefore, the Hamiltonian vectorfields of the functions
$q^{Z_i},h_1,h_2,k,f^{X_i}$ at the point $((v,z),V+Z)$ are of the
form
$$
\begin{gathered}
(Z_i,0),\\
(*\,,2\<V,E_1(Z)\>j(Z)E_1(Z)+2\<V,E_2(Z)\>j(Z)E_2(Z)),\\
(*\,,2\<V,E_3(Z)\>j(Z)E_3(Z)+2\<V,E_4(Z)\>j(Z)E_4(Z)),\\
(*\,,j(Z)Y(Z))=(*\,,0),\\
(W,0)+\Phi(Z)\cos(2\pi\<X_i\,,v_\x-C(Z)V_\y\>)\cdot 2\pi
(-L_{(v,z)*}{}^t C(Z)X_i, -X_i - j(Z){}^t C(Z)X_i)\\ =
(W-\Phi(Z)\cos(\ldots)\cdot 2\pi L_{(v,z)*}{}^t C(Z)X_i,0),
\end{gathered}
$$
with some $W\in\z$, where the last equality follows from
$$-j(Z){}^t C(Z)X_i ={}^t(j(Z)\restr\x){}^t C(Z)X_i=X_i\,.
$$ The
formulas for $q^{Z_j},q^{Z_k}, f^{X_j}$ are analogous. Thus, the
second components of the Hamiltonian vectorfields of the functions
$q^W,k,f^X$ vanish. Since the functions $q^W,h_1,h_2,k$ depend
only on the second component of $((v,z),V+Z)$, it follows
immediately that the latter functions Poisson commute with the
$q^W,k,f^X$. It only remains to show that $\{h_1,h_2\}=0$ and
$\{f^{X_i},f^{X_j}\}=0$: The derivative of $h_2$ in direction of
$X_{h_1}$ vanishes because the second component of $X_{h_1}$ at
$((v,z),V+Z)$ is in $\spann\{E_1(Z),E_2(Z)\}$ and thus orthogonal
to $E_3(Z)$ and $E_4(Z)$; the derivative of $f^{X_j}$ in direction
of $X_{f^{X_i}}$ vanishes because ${}^t C(Z)X_i\in\y$ has
vanishing $\x$-component.
\end{proof}

\Section{Structure of submanifolds foliated by generic closed
orbits} \label{sec:struct}

\noindent In this section, we will describe the submanifolds of
the unit tangent bundles foliated by continuous families of closed
geodesics in the two manifolds from Example~\ref{exple}. We will
consider only families most of whose geodesics have velocity
vectors satisfying a certain genericity condition.
The result (Corollary~\ref{cor:struct}) will nicely reflect the
(non)integrability properties of the geodesic flows established
in Section~\ref{sec:int}.



\begin{notation}
\label{decomp} Let $(\Gamma\backslash N,g):=
(\Gamma(j)\backslash N(j),g(j))$
and $(\Gamma'\backslash N',g'):=(\Gamma(j')\backslash N(j'),g(j'))$
be the two manifolds from
Example~\ref{exple}. Let $\gamma:\R\to(N^{(\prime)},g^{(\prime)})$
be a geodesic. Recalling Notation~\ref{not:twostep}, write
$\dot\gamma(0)=(\gamma(0),V+Z)$ for some $V\in\v$, $Z=Z_c\in\z$.
If $|c|>|c_k|>0$ then we write $V=V_{c_k}+V_{|c|}+V_0$ where
$V_\lambda$ denotes the component of $V$ in the
$(-\lambda^2)$-eigenspace of $j^{(\prime)}(Z)^2$. In what follows,
we will restrict our attention to geodesics $\gamma$ with
``generic'' velocity fields; by this, we mean that the vectors
$Z=Z_c$ and $V$ satisfy the following genericity condition:
\begin{equation}
\label{gener} |c|>|c_k|>0\mathrm{\ and\ } V_{c_k}\ne0,
V_{|c|}\ne0, V_0\ne0.
\end{equation}
Note that by the geodesic equations~\ref{lem:geodeqs}, this
property is invariant under the geodesic flow; so all
$\dot\gamma(t)$ will satisfy the corresponding condition if
$\dot\gamma(0)$ does so. Moreover, note that the set of tangent
vectors satisfying this genericity condition is open and dense in
the tangent bundle $TN^{(\prime)}$.
\end{notation}

\begin{remark}
Let $\gamma:\R\to(N^{(\prime)},g^{(\prime)})$ be a geodesic, and
let $\tau>0$. Then $\gamma$ will descend to a $\tau$-periodic
geodesic in the quotient manifold $(\Gamma^{(\prime)}\backslash
N^{(\prime)},g^{(\prime)})$ if and only if
$$
a:=\gamma(\tau)\gamma(0)\inv\in\Gamma^{(\prime)}\mathrm{\ and\ }
\dot\gamma(\tau)=L_{a*}\dot\gamma(0).
$$
Thus, if $\dot\gamma(0)= (\gamma(0),V+Z)$, then a necessary
condition for $\gamma$ to descend to a $\tau$-periodic geodesic is
$\dot\gamma(\tau)=(\gamma(\tau),V+Z)$ with the same vector
$V+Z\in\n$. Assuming this condition and the genericity
condition~(\ref{gener}) for $V+Z$, we will in the following lemma
compute the translational element $a=\gamma(\tau)\gamma(0)\inv$ in
terms of $\tau$, $V+Z$, and $\gamma(0)$. We first supply some
notation concerning eigenvectors of $j'(Z)^2$, analogous to
Notation~\ref{not:ourj}.
\end{remark}

\begin{notrems}
\label{not:ourjp} Let $j':\z\to\so(\v)$ be the second of the two
maps from Example~\ref{exple}. For $Z\in\z$ with
$Z=Z_c=c_iZ_i+c_jZ_j+c_kZ_k$ let
\begin{gather*}
E'_1(Z):=X_i,\quad E'_2(Z):=X_j\\
E'_3(Z):=|c|(c_jY_i-c_iY_j), \quad
E'_4(Z):=c_k(c_iY_i+c_jY_j)-(c_i^2+c_j^2)Y_k
\end{gather*}
and
$$Y(Z):=Y_c=c_iY_i+c_jY_j+c_kY_k.
$$
Note that we then have
\begin{equation}\label{jpeq}
\begin{gathered}
j'(Z)E'_1(Z)=c_kE'_2(Z),\quad
j'(Z)E'_2(Z)=-c_kE'_1(Z),\\
j'(Z)E'_3(Z)=|c|E'_4(Z),\quad j'(Z)E'_4(Z)=-|c|E'_3(Z),\\
j'(Z)Y(Z)=0.
\end{gathered}
\end{equation}
\end{notrems}

\begin{lemma}
\label{lem:comp} Let $V+Z\in\v\oplus\z$ satisfy the genericity
condition~$\mathrm{(\ref{gener})}$. Let
$\gamma:\R\to(N^{(\prime)},g^{(\prime)})$ be a geodesic with
$\dot\gamma(0)=((v,z),V+Z)$, where $(v,z)=\gamma(0)$. Let
$\tau>0$, and assume that $\dot\gamma(\tau)=(\gamma(\tau),V+Z)$.
Then the translational element $a:=\gamma(\tau)\gamma(0)\inv$ is
equal to
\begin{equation}
\begin{split}
\label{trsl} \biggl(\tau V_0\,,&\;
\tau\bigl(1+\frac{|V_{\perp}|^2}{2|c|^2}\bigr)\cdot Z_c +
\tau\beta\bigl(\alpha_2-\frac{c_k}{c_i^2+c_j^2}(x_ic_i+x_jc_j)\bigr)
\cdot\bigl(-c_jZ_i+c_iZ_j\bigr)\\ {}+
\tau&\Bigl(-\frac{|V_{c_k}|^2}{2c_k|c|^2}+\beta
\bigl(\alpha_4-\frac1{c_i^2+c_j^2}(x_ic_j-x_jc_i)
\bigr)\Bigr)\cdot\bigl(c_k(c_iZ_i+c_jZ_j)-(c_i^2+c_j^2)Z_k\bigr)\biggr)
\end{split}
\end{equation}
in $(N,g)$, respectively to
\begin{equation}
\begin{split}
\label{trslp} \biggl(\tau V_0\,,&\;
\tau\bigl(1+\frac{|V_{\perp}|^2}{2|c|^2}\bigr)\cdot Z_c +
\tau\beta\bigl(-|c|\alpha'_3+y_k-\frac{c_k}{c_i^2+c_j^2}(y_ic_i+y_jc_j)\bigr)
\cdot\bigl(-c_jZ_i+c_iZ_j\bigr)\\ {}+
\tau&\Bigl(-\frac{|V_{c_k}|^2}{2c_k|c|^2}+
\beta\bigl(\alpha'_4-\frac1{c_i^2+c_j^2}(y_ic_j-y_jc_i)
\bigr)\Bigr)\cdot\bigl(c_k(c_iZ_i+c_jZ_j)-(c_i^2+c_j^2)Z_k\bigr)\biggr)
\end{split}
\end{equation}
in $(N',g')$, where, using Notation~$\mathrm{\ref{decomp}(ii)}$,
$\mathrm{\ref{not:ourj}}$, and~$\mathrm{\ref{not:ourjp}}$, we
write $Z=Z_c$\,, $V_0=\beta Y(Z)$,
$V_{\perp}:=V_{c_k}+V_{|c|}=:\sum_{m=1}^4\alpha^{(\prime)}_mE^{(\prime)}_m(Z)$,
and $v=x_iX_i+x_jX_j+y_iY_i+y_jY_j+y_kY_k$. Moreover, we have
$\tau c_k\in 2\pi\Z$ and $\tau|c|\in 2\pi\Z$.
\end{lemma}

\begin{proof}
By the geodesic equations~\ref{lem:geodeqs} and our assumption on
$\dot\gamma(\tau)$, we have $e^{\tau j(Z)}V=V$. By the genericity
condition on $V$, this implies here that $e^{\tau j(Z)}=\Id_\v$\,;
in particular, $\tau c_k$ and $\tau|c|$ are in $2\pi\Z$.

In addition, assume for the moment that $\gamma(0)$ equals
$e=(0,0)$, the neutral element of $N^{(\prime)}$. In this
situation, one sees either by using formulas
from~\cite{Eberlein:1994} or by direct integration using the
geodesic equations~\ref{lem:geodeqs} and our explicit knowledge of
the action of $j^{(\prime)}(Z)$ on the three different eigenspaces
 of~$j^{(\prime)}(Z)^2$:
\begin{equation*}
\gamma(\tau)=\bigl(\tau V_0\,,\tau Z+\tau[V_0,j(Z)\inv
V_\perp]+\frac12\tau[j(Z)\inv
V_{c_k}\,,V_{c_k}]+\frac12\tau[j(Z)\inv V_{|c|}\,,V_{|c|}]\bigr)
\end{equation*}
in $(N,g)$, and the analogous formula for $(N',g')$ with $j'(Z)$
instead of $j(Z)$ and $\liebr'$ instead of $\liebr$, where
$j^{(\prime)}(Z)\inv$ denotes the inverse of
$j^{(\prime)}(Z)\restr{\spann\{V_0\}^\perp}$\,. Now if
$\gamma(0)=(v,z)$ is arbitrary, then
$\bar\gamma:=L_{(v,z)}\inv\circ\gamma$ is a geodesic as just
considered, with $\dot{\bar\gamma}(0)=((0,0),V+Z)$. Then
$\gamma(\tau)\gamma(0)\inv=(v,z)\bar\gamma(\tau)(v,z)\inv$. For
any element $(\bar v,\bar z)\in N^{(\prime)}$, we have $(v,z)(\bar
v,\bar z)(v,z)\inv=(\bar v,\bar z+[v,\bar v]^{(\prime)})$ by
\ref{not:twostep}(i). Thus, by adding the term $[v,\tau V_0]$ to
the $z$-component in the above formula, we get
\begin{equation*}
a=\bigl(\tau V_0\,,\tau Z+\tau[v,V_0]+\tau[V_0,j(Z)\inv
V_\perp]+\frac12\tau[j(Z)\inv
V_{c_k}\,,V_{c_k}]+\frac12\tau[j(Z)\inv V_{|c|}\,,V_{|c|}]\bigr)
\end{equation*}
in $(N,g)$, and the analogous formula for $(N',g')$.

The rest of the proof consists in evaluating this formula in
$(N,g)$ and $(N',g')$, respectively, using the definition of the
Lie brackets $\liebr$ and $\liebr'$ and the facts $V_0=\beta
Y(Z)=\beta Y_c$\,, $V_\perp=V_{c_k}+V_{|c|}$,
$V_{c_k}=\alpha^{(\prime)}_1E^{(\prime)}_1(Z)+\alpha^{(\prime)}_2E^{(\prime)}_2(Z)$,
$V_{|c|}=\alpha^{(\prime)}_3E^{(\prime)}_3(Z)+\alpha^{(\prime)}_4E^{(\prime)}_4(Z)$,
$j^{(\prime)}(Z)\inv
V_{c_k}=-\frac{\alpha^{(\prime)}_1}{c_k}E^{(\prime)}_2(Z) +
\frac{\alpha^{(\prime)}_2}{c_k}E^{(\prime)}_1(Z)$ and similarly
for $j^{(\prime)}(Z)\inv V_{|c|}$\,; for developing the resulting
$z$-component into the claimed form it is moreover useful to note
that
$Z_k=\frac{c_k}{|c|^2}Z_c-\frac1{|c|^2}\bigl(c_k(c_iZ_i+c_jZ_j)-
(c_i^2+c_j^2)Z_k\bigr)$.
The computation is a little tedious, but straightforward;
we spare the reader the details here.
\end{proof}

\begin{remark}
Note that by Lemma~\ref{lem:comp}, all translational elements
belonging to closed geodesics with velocity fields satisfying the
genericity condition~(\ref{gener}) are elements of the codimension
two submanifold
$\exp(\y\oplus\z)\subset\exp(\n^{(\prime)})=N^{(\prime)}$. So,
only a quite special type of free homotopy classes in
$\Gamma^{(\prime)}\backslash N^{(\prime)}$ contains closed
geodesics with generic velocity fields in the above sense. This,
however, should not lead to doubts as to whether the notion of
genericity is out of place here. The set of free homotopy classes
is discrete anyway, and there is no notion of genericity within
this set. Our genericity condition concerns only the velocity
vectors of the closed geodesics, as elements of the manifold
$TN^{(\prime)}$. It is a common phenomenon in compact Riemannian
nilmanifolds that ``generic'' closed geodesics belong only to a
quite special set of free homotopy classes. For example, in a
compact Riemannian Heisenberg manifold (or, more generally, in any
nonsingular two-step Riemannian nilmanifold), all ``generic'' closed
geodesics belong to central free homotopy classes, while the other
free homotopy classes contain only geodesics with very special
velocity fields.
\end{remark}

\begin{corollary}
\label{cor:struct} Under the assumptions of
Lemma~$\mathrm{\ref{lem:comp}}$, assume that
$a\in\Gamma^{(\prime)}$; in particular, $\gamma\restr{[0,\tau]}$
descends to a closed geodesic $\hat\gamma$ in
$(\Gamma^{(\prime)}\backslash N^{(\prime)}, g^{(\prime)})$.
Moreover, assume that $\gamma$ is a unit speed geodesic
$\mathrm($i.e., $|V+Z|=1\mathrm)$. Then the velocity fields of the
largest continuous family of closed unit speed geodesics
containing $\hat\gamma$ foliate a submanifold $C$, resp.~$C'$, of
the unit tangent bundle $S(\Gamma\backslash N)$, resp.
$S(\Gamma'\backslash N')$, with the following properties:
\begin{itemize}
\item[(i)] $C$ is diffeomorphic to
$T^6\times S^3$; it consists of a one-parameter family of
submanifolds invariant under the geodesic flow and diffeomorphic
to $T^6\times S^1\times S^1= T^8$, and two singular
seven-dimensional fibers diffeomorphic to $T^6\times S^1= T^7$.
The invariant $T^8$~fibers are level sets of the first integrals
from~Lemma~$\mathrm{\ref{lem:firstints}}$. The above decomposition
of~$C$ arises from the decomposition of the $S^3$ factor into a
one-parameter family of $T^2=S^1_r\times S^1_{\sqrt{1-r^2}}$
fibers and two singular $S^1$ fibers.
\item[(ii)] $C'$ is
diffeomorphic to $H^3\times T^3\times S^3$, where $H^3$ is a
compact three-dimensional Heisenberg manifold $\mathrm($in particular,
not diffeomorphic to~$T^3\mathrm)$. $C'$ consists of a
one-parameter family of submanifolds invariant under the geodesic
flow and diffeomorphic to $H^3\times T^3\times T^2=H^3\times T^5$,
and two singular seven-dimensional fibers diffeomorphic to
$H^3\times T^4$. This decomposition of~$C'$ arises from the
corresponding decomposition of the $S^3$ factor.
\end{itemize}
\end{corollary}

\begin{proof}
Let $\hat\gamma_s$ be a continuous family of closed unit speed
geodesics with $\hat\gamma_0=\hat\gamma$. For the time being, we
assume that the velocity field of each $\hat\gamma_s$ satisfies
the genericity condition~(\ref{gener}). Each $\hat\gamma_s$ has
the same length~$\tau$ as $\hat\gamma$ by the first variation
formula. We lift $\hat\gamma_s$ to a continuous family
$\gamma_s:[0,\tau]\to(N^{(\prime)},g^{(\prime)})$ of unit speed
geodesics in the universal cover such that
$\gamma_0(0)=\gamma(0)$, and extend them to geodesics
$\gamma_s:\R\to(N^{(\prime)},g^{(\prime)})$; thus
$\gamma_0=\gamma$. Since the family is continuous, the
translational element
$\gamma_s(\tau)\gamma_s(0)\inv\in\Gamma^{(\prime)}$ must be
constant in~$s$, hence equal to~$a$ for each~$s$. Writing
$\dot\gamma_s(0)=((v^s,z^s),V^s+Z^s)$ and $Z^s=Z_{c(s)}$\,, we
immediately read off from the first component of~$a$ in
Lemma~\ref{lem:comp} that $V_0^s\equiv V_0$, whence all $Y_{c(s)}$
are scalar multiples of each other. Moreover, by the lemma,
$\tau|c(s)|\in 2\pi\Z$ is constant. Therefore, $c(s)\equiv c$. By
$V_0=\beta^s Y_{c(s)}=\beta^sY_c$ we obtain $\beta^s\equiv\beta$.
From the coefficient of $Z_c$ occurring in the second component
of~$a$ in Lemma~\ref{lem:comp}, we read off that
$|V^s_\perp|\equiv|V_\perp|$. We continue the discussion
separately for the two parts of the statement:

(i) We obtain no restriction at all for the coefficients
$y^s_i,y_j^s,y_k^s,z_i^s,z_j^s,z_k^s$ of $v^s$ and $z^s$. There is
also no further restriction on $V^s_{c_k}$ and $V^s_{|c|}$ apart
from $|V^s_{c_k}|^2+|V^s_{|c|}|^2=|V_\perp|^2$ and (by the
genericity condition) the requirement that both vectors be
nonzero. These two vectors then determine the coefficients
$\alpha_2^s$ and $\alpha_4^s$\,. These, in turn, are by
Lemma~\ref{lem:comp} and the constancy of~$a$ seen to determine
$x_i^sc_i+x_j^sc_j$ and $x_i^sc_j-x_j^sc_i$\,, and thus to
determine $x_i^s$ and $x_j^s$ (recall that $(c_i,c_j)\ne0$).

Note that the sets of points in $\Gamma\cdot(\bar v,\bar z)\in
\Gamma\backslash N$ with fixed $x_i, x_j$ coordinates for~$\bar v$
and arbitrary $y_i,\ldots,z_k$ coordinates are diffeomorphic to
$T^6$ because the normal subgroup of~$N$ whose Lie algebra is
spanned by $\{Y_i,\ldots,Z_k\}$ is abelian. Moreover, the pairs of
possible vectors $V^s_{c_k}$ and $V^s_{|c|}$ fill out the sphere
$S^3$ with radius $|V_\perp|$ in $\spann\{V_0\}^\perp\subset\v$,
except for two $S^1$ fibers where the first or the second
component is zero. If we now drop the genericity condition (i.e.,
allow $V^s_{c_k}$ or $V^s_{|c|}$ to vanish), the possible pairs of
vectors fill out the entire~$S^3$.

Conversely, it is clear that any initial velocity $((v^s,z^s
),V^s+Z^s)\in TN$ obtained by choosing the various coordinates
according to the restrictions and degrees of freedom described
above will indeed correspond to a closed unit speed geodesic in
the maximal continuous family containing $\hat\gamma$. By the
above discussion, this shows that the corresponding
submanifold~$C$ of the unit tangent bundle is indeed diffeomorphic
to $T^6\times S^3$. (Note that quotienting by $\Gamma$ will not
identify pairs of tangent vectors with different
$(x_i,x_j)$-coordinates of the basepoint because such have been
seen to arise only from different velocity vectors.)

The decomposition into the two singular fibers $T^6\times S^1$ and
the one-parameter family of fibers $T^6\times S^1\times S^1$ is
respected by the geodesic flow because choosing one of these
fibers corresponds to fixing the norm of $|V^s_{c_k}|$ (and hence
of $|V^s_{|c|}|$); note that these norms are invariant under the
geodesic flow by~\ref{lem:geodeqs}.

Finally, the $T^6\times S^1\times S^1=T^8$ fibers turn out to be
level sets of the first eight integrals given in
Lemma~\ref{lem:firstints}: The coefficients $c_i,c_j,c_k$ of $Z$
which are constant here are the values of $q^{Z_i}, q^{Z_j},
q^{Z_k}$. The values $|V^s_{c_k}|^2$ and $|V^s_{|c|}|^2$, which
are constant in such a $T^8$ fiber, are, up to some multiplicative
constants depending on~$c$, just the values of $h_1$ and $h_2$.
The value of the first integral~$k$ is just $\<V^s,Y_c\>=\<\beta
Y_c\,,Y_c\>=\beta|c|^2$ which is constant as well. Moreover, one
straightforwardly computes
\begin{equation*}
\begin{gathered}
v_\x-C(Z)V_\y=x_iX_i+x_jX_j-\frac{\alpha_2}{c_k}E_1(Z)-\frac{\alpha_4}{|c|}E_3(Z)\\
=
\bigl(-\frac{\alpha_2}{c_k}+\frac1{c_i^2+c_j^2}(x_ic_i+x_jc_j)\bigr)E_1(Z)
+
\bigl(-\frac{\alpha_4}{|c|}+\frac1{|c|(c_i^2+c_j^2)}(x_ic_j-x_jc_i)\bigr)E_3(Z).
\end{gathered}
\end{equation*}
Comparing with the coefficients in the $z$-component of the
translational element~$a$ in Lemma~\ref{lem:comp}, we see that
constancy of those coefficients is equivalent to constancy of the
first integrals $f^{X_i}$ and $f^{X_j}$ (note that
$\Phi(Z)=\Phi(Z_c)$ is nonzero here because of $c_k\ne0$).
Conversely, one cannot continuously move out of such a $T^8$~fiber
without changing the values of any of the first integrals.
This finishes the proof of statement~(i).

(ii) The discussion is analogous to the corresponding discussion
in~(i). This time, we obtain no restriction for the coefficients
$x_i^s,x_j^s,z_i^s,z_j^s,z_k^s$ of $v^s$ and $z^s$\,, and the same
restrictions vs.~freedoms as in (i) for the two vectors
$V^s_{c_k}$ and $V^s_{|c|}$\,. The vector $V^s_{|c|}$ then
determine the coefficients $\alpha^{\prime\,s}_3$ and
$\alpha^{\prime\,s}_4$\,, which, in turn, determine
$y_k^s-\frac{c_k}{c_i^2+c_j^2} (y_i^sc_i+y_j^sc_j)$ and
$y_i^sc_j-y_j^sc_i$\,. The corresponding triples
$(y_i^s,y_j^s,y_k^s)$ constitute an affine line in the
$\y$-component of $\v=\x\oplus\y$ whose direction is rational with
respect to the sublattice of $\Gamma'$ spanned by $Y_i,Y_j,Y_k$\,:
Note that $c$ must be a scalar multiple of a rational vector
because $\tau V_0=\tau\beta Y_c$ was by Lemma~\ref{lem:comp}
the first component of $a\in\Gamma'$.

The sets of points in $\Gamma'\cdot(\bar v,\bar z)\in
\Gamma'\backslash N'$ with fixed $y_i, y_j, y_k$ coordinates for
$\bar v$ and arbitrary $x_i,x_j,z_i,z_j,z_k$ coordinates are
diffeomorphic to $H^3\times T^2$, where $H^3$ is a compact
three-dimensional Heisenberg manifold, corresponding to the
three-dimensional Heisenberg algebra spanned by $\{X_i,X_j,Z_k\}$
in~$\n'$. Moreover, the pairs of possible vectors $V^s_{c_k}$ and
$V^s_{|c|}$ have exactly the same freedom as in~(i) and fill out,
after dropping the genericity condition on these, the sphere $S^3$
of radius~$|V_\perp|$ in $\spann\{V_0\}^\perp$. Finally, the
affine line of possible $(y_i^s,y_j^s,y_k^s)$ determined via
$\alpha^{\prime\,s}_3$ and $\alpha^{\prime\,s}_4$ by these vectors
yields another $S^1$ factor after quotienting by~$\Gamma'$ (recall
that the line had rational direction); this indeed splits off as a
factor because $\y$ commutes with $\x\oplus\z\subset\n'$.

By analogous arguments as in (i), we conclude that $C'$ is indeed
diffeomorphic to $H^3\times S^1\times T^2\times S^3= H^3\times
T^3\times S^3$, decomposed as claimed in the statement.
\end{proof}

\begin{remark}
We do not present here the corresponding analysis for initial
velocities which do not satisfy the genericity condition~(\ref{gener}).
We just mention that the submanifolds of the unit tangent
bundle fibered by closed geodesics belonging to central free
homotopy classes but having noncentral velocity fields (this
is the case for most of such geodesics) are of dimension ten
for both of the manifolds from Example~\ref{exple}.
There are other special free homotopy classes for which
dimensions lower than eight occur.
\end{remark}

\begin{corollary}
\label{cor:dense}
For both of the $8$-dimensional manifolds from
Example~$\mathrm{\ref{exple}}$, the set of initial velocity
vectors of closed unit speed geodesics satisfying the genericity
condition~$\mathrm{(\ref{gener})}$ is dense in the unit tangent
bundle. In particular, by Corollary~$\mathrm{\ref{cor:struct}}$:
The $T^8$ fibers of the unit tangent bundle
of the first manifold $(\Gamma\backslash N,g)$
which are invariant under the geodesic flow and are
level sets of the eight first integrals from
Lemma~$\mathrm{\ref{lem:firstints}}$ are dense in the unit tangent
bundle; the same holds for the corresponding
invariant $H^3\times T^5$ fibers
in the unit tangent bundle of the second manifold
$(\Gamma'\backslash N',g')$.
\end{corollary}

\begin{proof}
We carry out the proof for $(\Gamma\backslash N,g)$ using
formula~(\ref{trsl}) from Lemma~\ref{lem:comp};
the proof for $(\Gamma'\backslash N',g')$
is completely analogous, using formula~(\ref{trslp}) instead.

Let $((\bar v,\bar z),\bar V+\bar Z)$ be an arbitrary element in
the unit tangent bundle of $(N,g)$. Write $\bar Z=Z_{\bar c}$ for some
$\bar c\in\R^3$. Let $\epsilon>0$ be arbitrary.
Choose $c\in\Q^3$ such that $|c-\bar c|<\epsilon$,
$|c|>|c_k|>0$, and $c_k/|c|\in\Q$. For the latter condition,
note that the rational points $c$ with rational value of
$c_k/|c|$ are dense in~$\R^3$ because the rational points
are dense in the unit sphere (since the standard stereographic
projection preserves rationality). Choose $\sigma>0$ such that
$\sigma c_k$ and $\sigma|c|$ are in $2\pi\Z$ (for instance,
$\sigma:=2\pi\frac q{|c|}$ if $c_k/|c|=p/q$ with $p\in\Z$, $q\in\N$).
Let $Z:=Z_c$ and write $\bar V=\bar V_0+\bar V_\perp$ with
$\bar V_0\in\spann\{Y_c\}$ and $\bar V_\perp \perp Y_c$\,.
Choose $0\ne V_0\in\spann\{Y_c\}$ such that $|V_0-\bar V_0|<\epsilon$
and $V_0\in\frac1\sigma\Q Y_c$\,. Choose $V_\perp\perp
Y_c$ such that $|V_\perp-\bar V_\perp|<\epsilon$, $V_\perp$ is
not orthogonal to neither $\spann\{E_1(Z),E_2(Z)\}$ nor
$\spann\{E_3(Z),E_4(Z)\}$, and such that $|V_\perp|^2\in
\frac1\sigma\Q$. Let $V:=V_0+V_\perp$\,. Finally, choose
$v\in\v$ such that $|v-\bar v|<\epsilon$ and
\begin{equation*}
\beta\bigl(\alpha_2-\frac{c_k}{c_i^2+c_j^2}(x_ic_i+x_jc_j)\bigr)\in
\frac1\sigma\Q\mathrm{\ \ and\ \ }
-\frac{|V_{c_k}|^2}{2c_k|c|^2}+\beta\bigl(\alpha_4-
\frac1{c_i^2+c_j^2}(x_ic_j-x_jc_i)\bigr)\in \frac1\sigma\Q,
\end{equation*}
where $V_0=\beta Y_c$\,, $V_\perp=\sum_{m=1}^4\alpha_mE_m(Z)$,
and $v=x_iX_i+x_jX_j+y_iY_i+y_jY_j+y_kY_k$\,. Let $\bar z:=z$.
Then by Lemma~\ref{lem:comp}, for the geodesic $\gamma$ with initial
velocity $((v,z),V+Z)$ we have $\gamma(\sigma)\gamma(0)\inv\in
\exp(\spann_\Q\{Y_i,\ldots,Z_k\})$. By replacing $\sigma$ with
a suitable multiple $\tau=m\sigma$, where $m\in\N$, we
obtain $a:=\gamma(\tau)\gamma(0)\inv\in\Gamma$ and (still)
$\tau c_k\in 2\pi\Z$, $\tau|c|\in 2\pi\Z$. In particular,
$e^{\tau j(Z)}=\Id_\v$ and hence $L_{a*}\dot\gamma(0)=\dot\gamma(\tau)$,
thus $\gamma\restr{[0,\tau]}$ descends to a $\tau$-periodic
geodesic in $(\Gamma\backslash N,g)$.
Summarizing, we have shown: Arbitrarily close to any unit tangent
vector $((\bar v,\bar z),\bar V+\bar Z)$, we find another tangent
vector $((v,z),V+Z)$ satisfying the genericity condition~(\ref{gener})
corresponding to the initial vector of a closed geodesic in
$(\Gamma\backslash N,g)$. Choosing a sequence of such vectors
converging to $((\bar v,\bar z),\bar V+\bar Z)$ and normalizing
each vector of this sequence, we obtain a sequence of unit vectors
converging to our given unit vector, consisting of initial
velocities of closed geodesics and satisfying the genericity
condition~(\ref{gener}).
This proves the statement.
\end{proof}

\begin{remark}
\label{rem:CIH}
On the one hand, it seems surprising at first sight
that two manifolds with the same Laplace
spectrum, thus with the same wave trace, can
differ so radically with respect to the behaviour of their
geodesic flows. On the other hand, as explained in the Introduction,
the fact that continuous families of closed geodesics in our
pair of manifolds have been seen to fill out larger dimensional
families than Lagrangian tori shows that the suitable version of the
condition of ``clean simple length spectrum'', under which
the problem of audibility of complete integrability and any positive
conjecture in this respect would make most sense, is violated here anyway.
This yields an explanation of why isospectrality in these
examples fails to entail similar integrability properties.

An even more obvious explanation would be provided if the
two manifolds failed to satisfy the so-called Clean Intersection
Hypothesis, a condition always needed for even establishing
a wave trace formula; see~\cite{DG:1975}. However, it turns out that both of
our manifolds actually do satisfy the Clean Intersection
Hypothesis, which is defined as the following condition:
For any number $\ell>0$
and any free homotopy class $\alpha$ containing a closed geodesic
of length~$\ell$, the subset $W_\ell(\alpha)$ of the unit tangent
bundle consisting of the velocity vectors of all unit speed closed
geodesics of length~$\ell$ and contained in~$\alpha$ is a finite
union of submanifolds of the unit tangent bundle, and for the
differential of the time-$\ell$-map of the geodesic flow at any
point of $W_\ell(\alpha)$, the
eigenspace associated with the eigenvalue~$1$
is not larger than the tangent space of $W_\ell(\alpha)$
at that point.

In fact, Ruth Gornet~\cite{Gornet:2005} has recently shown that a
compact, two-step Riemannian nilmanifold $(\Gamma\backslash N,g)$
satisfies this hypothesis if and only if, using the notation
from~\ref{def:twostep} and~\ref{notrems:latt}, for all
$V+Z\in\exp\inv(\Gamma)\subset\n$ and all nonzero eigenvalues~$\pm
i\theta$ of $j(\proj_{[V,\n]^\perp}Z)$ (where $\proj$ denotes
orthogonal projection) we have $\theta\notin\pi\Q$. In our
manifolds from Example~\ref{exple}, let us endow $\v$ and $\z$ with
the rational structure given by $\spann_{\Q}\{X_i,\ldots,Y_k\}$
and $\spann_\Q\{Z_i,Z_j,Z_k\}$. Then for all
$V+Z\in\exp\inv(\Gamma^{(\prime)})$, $V$ and $Z$ are rational
vectors by our definition of~$\Gamma^{(\prime)}$. The subspace
$[V,\n]=[V,\v]\subseteq\z$ is then rational as well. Elementary
arguments (using the orthonormality of the rational basis
$\{Z_i,Z_j,Z_k\}$)
show that the orthogonal projection of the rational
vector~$Z$ to this rational subspace is again a rational vector.
Denote the resulting vector by $Z_V^\perp=c_iZ_i+c_jZ_j+c_kZ_k$.
Then $c_i,c_j,c_k$ are rational; in particular, the possibly
nonzero eigenvalues $-c_k^2$ and $-|c|^2$ of $j(Z_V^\perp)^2$ are
rational. So these can, if nonzero, never be in $\pi^2\Q$. This
implies that Gornet's necessary and sufficient criterion for the
Clean Intersection Hypothesis is indeed satisfied.
\end{remark}

\end{document}